\documentclass[12pt]{iopart}
\usepackage{iopams, graphicx, cite, color}%
\newtheorem{thm}{Theorem}[section]

\newtheorem{lem}[thm]{Lemma}

\newcommand{\Real}{\mathbb{R}}
\newcommand{\eps}{\varepsilon}
\newcommand{\dmn}{\mathop{\rm dom}}
\newcommand{\supp}{\mathop{\rm supp}}
\newcommand{\dom}{\mathop{\rm dom}}
\newcommand{\norm}[1]{\left\Vert#1\right\Vert}
\newcommand{\tfrac}[2]{\textstyle{\frac{#1}{#2}}}
\newcommand{\abs}[1]{\left\vert#1\right\vert}

\newcommand{\bl}[1]{{\color{black}{#1}}}
\newcommand{\rd}[1]{{\color{black}{#1}}}
\renewcommand{\emph}[1]{{\textit{#1}}}

\begin{document}
\title[On norm resolvent convergence of  Schr\"{o}dinger operators]{On norm resolvent convergence of  Schr\"{o}dinger operators with $\delta'$-like potentials}
\author{Yu D Golovaty$^1$ and R O Hryniv$^{1,2,3}$}
\address{$^1$ Department of Mechanics and Mathematics, Ivan Franko National University of Lviv, 1 Universytetska str., 79000 Lviv, Ukraine }

\address{$^2$ Institute for Applied Problems of Mechanics and Mathematics, 3b Naukova str.,
79601 Lviv, Ukraine}

\address{$^3$ Institute of Mathematics, the University of Rzesz\'{o}w, 16\,A Rejtana str., 35-959 Rzesz\'{o}w, Poland}

\ead{yu\_\,holovaty@franko.lviv.ua and rhryniv@iapmm.lviv.ua}

\begin{abstract}
For a function~$V\,:\,\Real \to \Real$ that is integrable and compactly supported, we prove the norm resolvent convergence, as $\eps\to0$, of a family $S_\eps$ of one-dimensional Schr\"odinger operators on the line of the form
\[
    S_\eps:= -\frac{\rmd^2}{\rmd x^2} + \frac1{\eps^2}V\Bigl(\frac{x}{\eps}\Bigr).
\]
If the potential~$V$ satisfies the conditions
\[
 \int_\Real V(\xi)\,\rm d\xi=0,\qquad \int_\Real \xi V(\xi)\,\rm d\xi=-1,
\]
then the functions $\eps^{-2}V(x/\eps)$ converge in the sense of distributions as $\eps\to 0$ to~$\delta'(x)$, and the limit $S_0$ of~$S_\eps$ might be considered as a `physically motivated' interpretation of the one-dimensional Schr\"odinger operator with potential $\delta'$. In 1985, \v{S}eba  claimed that the limit operator $S_0$ is the direct sum of the free Schr\"odinger operators on positive and negative semi-axes subject to the Dirichlet condition at~$x=0$, which suggested that in dimension~$1$ there is no non-trivial Hamiltonian with potential $\delta'$. In this paper, we show that in fact~$S_0$ essentially depends on~$V$: although the above results are true generically, in the exceptional (or `resonant') case, the limit~$S_0$ is non-trivial and is determined by the properties of an auxiliary Sturm--Liouville spectral problem associated with~$V$. We then set $V(\xi) = \alpha \Psi(\xi)$ with a fixed $\Psi$ and show that there exists a countable set of  resonances $\{\alpha_k\}_{k=-\infty}^\infty$ for which a partial transmission of the wave package occurs for $S_0$.
\end{abstract}



\section{Introduction}

Assume that~$V$ is an integrable function of compact support contained in the interval~$[-1,1]$ and consider Schr\"odinger operators $S_\eps$ on the real line given by
\begin{equation}\label{eq:intr.Seps}
    S_\eps = -\frac{\rmd^2}{\rmd x^2} + \frac1{\eps^2} V\Bigl(\frac{x}{\eps}\Bigr).
\end{equation}
Here~$\eps$ is a positive parameter, and one of the questions of our primary interest in this paper is the behaviour of~$S_\eps$ as~$\eps$ tends to zero.

The motivation for this question stems from the fact that if $V$ has zero mean and its first moment $\int_\Real \xi V(\xi)\,\rmd \xi$ is $-1$, then the functions~$\eps^{-2}V(x/\eps)$ converge in the sense of distributions as~$\eps\to0$ to the derivative~$\delta'$ of the Dirac delta-function. Therefore, if the Hamiltonians~$S_\eps$ converge (in some topology) as~$\eps\to0$ to a limit $S_0$, then it is natural to regard~$S_0$ as a realization of the Schr\"odinger operator with a potential~$\delta'$.

Schr\"odinger operators in~$\Real^d$, $d\geq1$, with singular distributional potentials supported on a discrete set (such potentials are usually termed ``point interactions'') have attracted considerable attention both in the physical and mathematical literature over several past decades. One of the reasons for this is that such singular Hamiltonians have widely been used in quantum mechanics to model interactions in particles and atoms. For instance, as back as in~1931 Kronig and Penney~\cite{KronigPenney} already suggested their model of electrons moving in the crystal lattices that used one-dimensional Schr\"odinger operators with periodic potentials of the form $\alpha\sum_{n\in\mathbb Z}\delta(\cdot-n)$ as the corresponding Hamiltonians. Another reason is that Schr\"odinger operators with point interactions often form ``solvable'' models in the sense that the resolvents and other objects related to such operators can explicitly be calculated (see the books by Albeverio, Gesztesy, H{\o}egh-Krohn, and Holden~\cite{Albeverio2edition} and by Albeverio and Kurasov~\cite{AlbeverioKurasov} discussing point interactions and more general singular perturbations of the free Schr\"odinger operators in~$\Real^d$ and the extensive bibliography lists therein). This in turn allows one to discover new, unusual effects not seen for regular potentials as was the case, e.g.\ with the singular Wannier--Stark systems~\cite{Exner:1995}.

One of the most natural ways to define a Hamiltonian corresponding to a point interaction supported by $x=0$ is to first restrict the free Schr\"odinger operator onto the set of smooth functions that vanish in a neighbourhood of $x=0$ and then take a self-adjoint extension of the resulting symmetric operator~$S_0'$. This approach was suggested for the first time by Berezin and Faddeev~\cite{BerezinFadeev} in 1961 for three-dimensional Schr\"odinger operators with a potential~$\alpha\delta$; the operator~$S_0'$ then has deficiency indices~$(1,1)$, and all self-adjoint extensions of $S_0'$ form a one-parameter family $H_\gamma$. In~\cite{BerezinFadeev}, the authors suggested a physically motivated choice of $\gamma$ for a given~$\alpha$ based on the renormalization technique pertinent to the physical setting of the problem. The existence of nontrivial limits of Hamiltonians as the singular potential $\alpha\delta$ is approximated by a regular sequence of compactly supported ones was also demonstrated in~\cite{Friedman:1972}. We also notice that differential operators with singular coefficients can be studied within the framework of the contemporary theory of new generalized functions such as the Colombeau or Egorov's algebras, in which multiplication is well defined and which contain rich sets of (different) $\delta$-functions obtained via approximating sequences with different profiles. This approach was applied in \cite{Antonevich:1999} to study the Schr\"odinger operator in $\Real^3$ with $\delta$-potentials.

In dimension~$1$, the symmetric operator~$S_0'$ has deficiency indices~$(2,2)$; thus the set of all its self-adjoint extensions forms a four-parameter family, and the problem of choosing a single Hamiltonian corresponding to a particular singular point interaction becomes more subtle. There are two families of extensions that have been studied especially well and have a clear physical interpretation.

The first one, denoted by $S_{\alpha,\delta}$, is given by $S_{\alpha,\delta}f=-f''$ on the domain
\[
\fl
    \dom S_{\alpha,\delta} := \{f \in W^{2}_2(\Real\setminus\{0\}) \mid
                f(0+) = f(0-)=: f(0), \quad
                f'(0+) - f'(0-)= \alpha f(0) \}
\]
and corresponds to the Hamiltonian
\[
     - \frac{\rmd^2}{\rmd x^2} + \alpha \delta(x).
\]
That such an identification is most natural is seen from the fact that various known approaches to definition of the above Hamiltonian (e.g.\ via the form sum~\cite{AlbeverioKoshmanenko,Koshmanenko:1999}, generalized sum~\cite{BrascheNizhnik:2002,AlbeverioNizhnik:2000}, approximation by regular potentials~\cite{BrascheFigariTeta}, and regularization by quasi-derivatives~\cite{ShkalikovSavchukTMMO2003}) lead to the same operator~$S_{\alpha,\delta}$. In particular, if~$V_n$ are regular functions converging weakly to the measure~$\alpha\delta$, then the corresponding regular Schr\"odinger operators
\begin{equation}\label{eq:Vn}
    - \frac{\rmd^2}{\rmd x^2} + V_n
\end{equation}
converge in the norm resolvent sense to~$S_{\alpha,\delta}$~\cite{BrascheFigariTeta}, \cite[Ch.~I.3.2]{Albeverio2edition}.

The second well-studied family, denoted $S_{\beta,\delta'}$, of self-adjoint extensions of the operator $S_0'$ in one dimension is given by~$S_{\beta,\delta'}f = -f''$ on the set of functions
\[\fl
    \dom S_{\beta,\delta'}= \{f \in W^{2}_2(\Real\setminus\{0\}) \mid
                f'(0+) = f'(0-)=: f'(0), \quad
                f(0+) - f(0-)= \beta f'(0) \}
\]
and is widely accepted as a model for Schr\"odinger operators with~$\delta'$-\emph{interactions} (not $\delta'$-potentials!), i.e.\ for the operators
\[
    - \frac{\rmd^2}{\rmd x^2} + \beta \langle\,\cdot\, , \delta'\rangle \delta'.
\]
Here $\langle f, \psi \rangle$ denotes the action $\psi(f)$ of a distribution~$\psi$ on a test function~$f$.
Physically, this model corresponds to an idealized dipole of zero range at $x=0$.
Similarly to the case of $\delta$-potentials, the operator~$S_{\beta,\delta'}$ can be given as the limit in the strong resolvent sense of regular rank-$1$ perturbations of the free Schr\"odinger operator by $\beta (\cdot, V_n)_{L_2} V_n$ whenever $V_n$ converge to $\delta'$ in the distributional sense~\cite{Albeverio2edition,AlbeverioKurasov}. Recently it was realized that $S_{\beta,\delta'}$ is the limit in the norm resolvent sense of the operators~\eref{eq:Vn}, under the special choice of~$V_n$~\cite{ExnerNeidhardtZagrebnov,AlbeverioNizhnik:2000}.

We note that one can also define Schr\"odinger operators with $\delta'$-interactions supported by subsets~$\Gamma$ of~$\mathbb R$. Spectral properties of the corresponding Hamiltonians for discrete~$\Gamma$ were studied e.g.\ in~\cite{GesHolJPA,KostenkoMalamud:2009} and for $\Gamma$ a Cantor set in~\cite{NizhFAA2003}. Also, there are papers where more general rank-$1$ perturbations of abstract positive Hamiltonians~$H_0$ are considered, namely perturbations by $\alpha (\cdot,\psi)\psi$, where $\psi$ is a generalized function in the domain of the operator $(H_0+I)^{-n}$, for $n\in\mathbb N$. The efficient way to treat such singular perturbations is via Krein's resolvent identity, see~\cite{AlbeverioKoshmanenko:1999a} and especially the papers by Kurasov~\cite{Kurasov:2003} and Nizhnik~\cite{NizhFAA2006} and the references therein.

The situation is more obscure with definition of the Schr\"odinger operator with a \emph{potential}~$\delta'$, which might be considered as an alternative way of modelling idealized dipole of zero range at~$x=0$. Mathematically, such a potential should be interpreted as an operator of multiplication by the distribution~$\delta'$. The product $\delta' f$ for smooth enough~$f$ is a distribution, so that any Schr\"odinger operator chosen to serve as a putative model for such a physical system should be understood in the distributional sense. A natural approach to defining the operator is to approximate $\delta'$  by regular potentials (e.g. in the distributional sense) and then to investigate the convergence of the corresponding family of regular Schr\"odinger operators. As we explained above, one possible way
to do this is to consider the family $S_\eps$ of~\eref{eq:intr.Seps} with $V$ of zero mean and nonzero first moment.

This approach was realized by \v{S}eba~\cite{SebRMP}. In fact, he considered a family of Schr\"odinger operators with more general short-range potentials~$\eps^{-\alpha} V(x/\eps)$ in theorem~4 of~\cite{SebRMP}
and claimed in particular that for $\alpha>\tfrac32$ the operators $S_\eps$ converge in the norm resolvent sense as $\eps\to0$ to the direct sum $S_-\oplus S_+$ of the unperturbed half-line Schr\"odinger operators subject to the Dirichlet boundary conditions at~$x=0$. One therefore could suggest that, in dimension one, no non-trivial interpretation of $\delta'$ potentials is possible. However, a careful analysis of the arguments given in~\cite{SebRMP} reveals that, although the proof is solid when $\alpha<2$, the term $\eps^{2-\alpha}T_\eps(k)$ in the expansion of the resolvent $(S_\eps - k^2)^{-1}$ does not vanish in the limit as $\eps \to0$ if~$\alpha\ge2$, thus putting under question the validity of the result in this case.

This was one of the motivations for us to re-examine the above convergence result.
In this paper, we confirm part of the statement from~\cite{SebRMP} that for an arbitrary real-valued integrable potential~$V$ of compact support the family~\eref{eq:intr.Seps} of Schr\"odinger operators~$S_\eps$ converges in the norm resolvent sense to a limit denoted~$S_0$. However, the limit~$S_0$ turns out to heavily depend on weather or not~$V$ is resonant. To define this notion, we consider the Sturm--Liouville operator~$\mathcal N$ given by
\begin{equation}\label{eq:N}
     {\mathcal N} y = -y'' + V y
\end{equation}
on the interval~$[-1,1]$ subject to the Neumann boundary conditions~$y'(-1) = y'(1)=0$ and call the potential~$V$ \emph{resonant} if $\mathcal N$ has a non-trivial null-space and \emph{non-resonant} otherwise. We show that in the non-resonant case,
the limiting operator~$S_0$ is just the direct sum $S_-\oplus S_+$ of the free Schr\"odinger operators on $\mathbb R_-$ and $\mathbb R_+$ respectively subject to the Dirichlet boundary condition at $x=0$, as claimed in~\cite{SebRMP}. In the resonant case, we take an eigenfunction~$u$ corresponding to the eigenvalue~$\lambda=0$, set%
\begin{footnote}
{Note that in the resonant case any solution of the equation $u''=Vu$ with $u'(-1)=0$ is constant outside the convex hull of the support of~$V$; thus the number~$\theta$ does not change if in the definition of the operator~$\mathcal N$ we replace the interval~$[-1,1]$ by any interval containing the support of~$V$.}
\end{footnote}%
$\theta:=u(1)/u(-1)$, and denote by $S(\theta)$ the free Schr\"odinger operator restricted to functions in $W_2^2(\Real\setminus\{0\})$ obeying the interface condition
\[
    f(0+) = \theta f(0-),
        \quad \theta f'(0+) = f'(0-);
\]
then the limit $S_0$ is given by~$S(\theta)$.

In fact, resonances in the transmission probability for $\delta'$-like potentials have earlier been observed by Christiansen~a.o.~\cite{ChristianZolotarIermak03}.
In that paper, an exactly solvable model~\eref{eq:intr.Seps} with a specially chosen step function~$V= \alpha \Psi$ was considered. The authors found a discrete set of intensities $\alpha_n$ (called resonant values) for which partial transmission through the limiting $\delta'$-potential occurs; this transmission was shown to rapidly decay as $\alpha_n$ becomes larger. The values $\alpha_n$ are roots of a transcendent equation depending on the regularization of $\delta'$. The findings of~\cite{ChristianZolotarIermak03} are also in contradiction with the results of~\cite{SebRMP} suggesting that the $\delta'$-barrier is completely opaque, no matter what value~$\alpha$ takes.
Exactly solvable models with other piecewise constant potentials as well as non-rectangular regularizations of $\delta'$ have later been studied in~\cite{ZolotarChristianIermak06,ZolotarChristianIermak07,ToyamaNogami}.
It is also worth mentioning the recent results by Zolotaryuk~\cite{Zolotaryuk08,Zolotaryuk09,Zolotaryuk10}.

In \cite{GolovatyManko1} a similar effect was discovered for the family of Schr\"odinger operators on the line of the form
\[
    -\frac{\rmd^2}{\rmd x^2} + \frac\alpha{\eps^2} \Psi\Bigl(\frac{x}{\eps}\Bigr)
        +W(x),
\]
where $\Psi\in C^\infty_0(-1,1)$, $W$ is a real valued potential tending to $+\infty$ as $|x|\to\infty$, and $\alpha\in\Real$ is a coupling constant. The map assigning a self-adjoint extension $S(\alpha,\Psi)$ of the operator $S_0'+W$ to each pair $(\alpha,\Psi)$ was constructed there. The choice of the extension is determined by proximity
of the energy levels and the pure states for the Hamiltonians with smooth and singular potentials respectively.
Two spectral characteristics of the profile $\Psi$ are introduced in \cite{GolovatyManko1}: the \emph{resonance set} $\Sigma_\Psi$, which is the $\alpha$-spectrum of the Sturm-Liouville problem $-w''+\alpha \Psi w=0$ on the interval $(-1,1)$ subject to the boundary conditions $w'(-1)= w'(1)=0$, and the \emph{coupling function} $\theta_\Psi
\colon\Sigma_\Psi \to \mathbb{R}$ defined via $\theta_\Psi (\alpha):=w_\alpha(1)/w_\alpha(-1)$, where
$w_\alpha$ is an eigenfunction corresponding to the eigenvalue $\alpha\in \Sigma_\Psi$.
In the case when the coupling constant~$\alpha$ does not belong to the resonance set, $S(\alpha,\Psi)$ is just the direct sum  of the  Schr\"odinger operators with a potential $W$ on semi-axes, subject to the Dirichlet boundary condition at the origin. In the resonant case, when $\alpha\in \Sigma_\Psi$,  $S(\alpha,\Psi)$ acts via $S(\alpha,\Psi)f = -f'' + Wf$ on an appropriate set of functions obeying the interface condition $f(0+)=\theta_\Psi (\alpha)f(0-)$ and $\theta_\Psi (\alpha)f'(0+)=f'(0-)$.

After we have established the main results of this paper, Prof.\ Albeverio drew our attention to the related work~\cite{AlbeverioCacciapuotiFinco:2007,CacciapuotiExner:2007,Seba:1985}. In~\cite{Seba:1985}, \v{S}eba
demonstrated existence of `resonant' non-trivial limits for a similar family of the Dirichlet Schr\"odinger operators on the half-line producing in the limit the Robin boundary condition at~$x=0$. When studying the problem of approximating a smooth quantum waveguide with a quantum graph, the authors of~\cite{AlbeverioCacciapuotiFinco:2007,CacciapuotiExner:2007} also faced the question on the norm resolvent convergence of the family $S_\eps$ of~\eref{eq:intr.Seps}. Under the assumptions that $V$ decays exponentially fast at~$\pm\infty$ and the mean value $\int_{\Real} V$ of $V$ is non-zero, the authors singled out the set of resonant potentials~$V$ producing a non-trivial limit of $S_\eps$ in the norm resolvent sense as $\eps \to0$. Although the papers~\cite{AlbeverioCacciapuotiFinco:2007,CacciapuotiExner:2007} treat a more general situation, the very important case $\int_{\Real} V=0$ giving $\delta'$-type potentials in the limit is excluded from the consideration, and the analysis there crucially relies on the fact that $\int_{\Real} V\ne 0$. In contrast, our approach is insensitive to the mean value of~$V$. We also note that our definition of the resonant potential agrees with that of the papers~\cite{AlbeverioCacciapuotiFinco:2007,CacciapuotiExner:2007}.

One should keep in mind that the results of the cited papers~\cite{AlbeverioCacciapuotiFinco:2007,CacciapuotiExner:2007,GolovatyManko1,ChristianZolotarIermak03} and of this paper do not allow to \emph{define} the Schr\"odinger operator with potential $\delta'$. Indeed, the limiting operator $S_0$ is shape-dependent, i.e.\ it depends on the background potential~$V$, so that different approximating families produce different limits. However, these results suggest the best choice of an idealized solvable quantum-mechanical model for a realistic quantum-mechanical device modelling a dipole.

We should also note that it is somehow surprising that the limit of~$S_\eps$ exists in the norm resolvent sense, although the potentials $V_\eps$ in general do not converge even in the distributional sense. Although $V_\eps(x)\to0$ for $x\ne0$ as $\eps\to0$, the topology of pointwise convergence is too weak to force any type of convergence of $S_\eps$.

The paper is organized as follows. In the next section, we consider the resonant case and prove the convergence result; in section~\ref{sec:non-res}, the easier non-resonant case is treated. In section~\ref{sec:delta'}, we specialize the above results to the case of $\delta'$-like potential~$V = \alpha \Psi$ and characterize the resonant values of~$\alpha$ (i.e., those $\alpha$ for which $V$ is resonant). Finally, in the last two sections, we discuss scattering at the resonant values of $\alpha$ and illustrate this effect by a simple example.

\emph{Notation.} Throughout the paper, $W_2^j(\Omega)$, $j=0,1,2$, stands for the Sobolev space of functions defined on a set $\Omega\subset\Real$ that belong to $L_2(\Omega)$
together with their derivatives up to order $j$. The norm in $W_2^2(\Omega)$ is given by
\[
    \|f\|_{W_2^2(\Omega)}
        := \bigl(\|f''\|^2_{L_2(\Omega)} + \|f\|^2_{L_2(\Omega)} \bigr)^{1/2},
\]
where
\[
    \|f\|_{L_2(\Omega)}:= \Bigl(\int_\Omega |f(x)|^2\,dx\Bigr)^{1/2}
\]
is the usual $L_2$-norm. We shall write $\|f\|$ instead of $\|f\|_{L_2(\Real)}$ and note that $\|f'\|_{L_2(\Omega)} \le \|f\|_{W_2^2(\Omega)}$ by interpolation.

\section{Resonant Case}\label{sec:resonance}

In this section, we analyze the more difficult resonant case where $\lambda=0$ is an eigenvalue of the Sturm--Liouville operator~$\mathcal{N}$ on~$(-1,1)$ defined in the introduction via~\eref{eq:N} and the Neumann boundary conditions at $x=\pm1$, and denote by $u$ a corresponding eigenfunction satisfying the condition $u(-1) =1$. Since the spectrum of~$\mathcal{N}$ is simple, the function $u$ is uniquely defined. Next, we set $\theta:=u(1)$ and denote by $S(\theta)$ the free Schr\"odinger operator on the line acting via $S(\theta)y=-y''$ on the domain
\begin{equation}\label{DomS(theta)}
    \dmn S(\theta) = \{ y \in W_2^2(\Real\setminus\{0\}) \mid y(0+) = \theta y(0-),
        \ \theta y'(0+) = y'(0-) \}.
\end{equation}
The operator~$S(\theta)$ is a symmetric extension by two dimensions of the symmetric operator~$S_0'$ with deficiency indices~$(2,2)$ and therefore is self-adjoint.

Denote by $v$ a solution of the Cauchy problem:
  \begin{equation}\label{EqV}
    -v''+V(\xi)v=0,\quad \xi\in(-1,1), \qquad v(-1)=0, \quad v'(-1)=1.
  \end{equation}
Clearly, we have $v\in W^2_2(-1,1)$; moreover, the Lagrange identity yields $(v' u - vu') \bigl|_{-1}^{1} =0$, i.e., $v'(1) = \theta^{-1}$.

\begin{lem}\label{LemmaResonant}
Fix $k^2 \in \mathbb{C}$ with $\Im k^2 \ne0$; then there exists $c>0$ with the property that for every
$f\in L_2(\Real)$ and $\eps>0$ there is $q_\eps\in L_2(\Real)$ with $\norm{q_\eps}\leq c\eps^{1/2}\norm{f}$ such that the function $y_\eps=(S(\theta)-k^2)^{-1}f+q_\eps$ belongs to $\dmn S_\eps$ and satisfies the inequality
\begin{equation}\label{SepsEst}
    \norm{(S_\eps-k^2)y_\eps-f}\leq c\eps^{1/2}\norm{f}.
\end{equation}
\end{lem}
\smallskip

\noindent
\textbf{Proof.} \rd{Fix an arbitrary $f\in L_2(\Real)$, set $y:=(S(\theta)-k^2)^{-1}f$ and
consider the auxiliary Cauchy problem}
\begin{eqnarray}\label{WepsProblem}
\fl \bl{ -w_\eps''+V(\xi)w_\eps=f(\eps\xi),\quad \xi\in(-1,1), \qquad w_\eps(-1)=0, \quad w_\eps'(-1)=0.}
\end{eqnarray}
We extend $v$, \bl{$w_\eps$} and the eigenfunction $u$ to the whole line by zero and introduce the function
\begin{equation*}
    z_\eps(x)=\bigl(1-\chi(x/\eps)\bigr)y(x)+y(0-)u(x/\eps)+\eps y'(0-)v(x/\eps)\bl{+\eps^2 w_\eps(x/\eps)},
\end{equation*}
with $\chi$ being the characteristic function of the interval~$[-1,1]$.
By construction, the function~$z_\eps$ belongs to $W_2^2(\Real\setminus\{-\eps,\eps\})$.  Although $z_\eps$ is in general discontinuous at the points $x=-\eps$ and $x=\eps$, its jumps and the jumps of its first derivative at these points are small.

To justify this, let $[g]_a$ denote the jump of a function~$g$ at a point $x=a$. We first observe that $(S(\theta)-k^2)^{-1}$ is a bounded operator from~$L_2(\Real)$ to the domain of~$S(\theta)$ equipped with the graph norm; since the latter space is equivalent to $W_2^2(\Real\setminus \{0\})$, there exists a constant $c_1>0$ independent of~$f$ such that
\begin{equation}\label{SinvEst}
    \|y\|_{W_2^2(\Real\setminus \{0\})}\leq c_1\|f\|.
\end{equation}
Now for $x=-\eps$ we get the estimates
\begin{eqnarray}\label{JumpZAt-Eps}
\fl \bigr|[z_\eps]_{-\eps}\bigl|
    =\bigr|y(0-)-y(-\eps)\bigl|\leq \int^0_{-\eps}|y'(t)|\,\rmd t
    \leq \eps^{1/2}\|y\|_{W_2^2(\Real\setminus \{0\})}
    \leq c_1\eps^{1/2}\|f\|,\\\label{JumpZ'At-Eps}
  \fl \bigr|[z'_\eps]_{-\eps}\bigl|=
    \bigr|y'(0-)-y'(-\eps)\bigl|
    \leq\int^0_{-\eps}|y''(t)|\,\rmd t
    \leq \eps^{1/2}\|y\|_{W_2^2(\Real\setminus \{0\})}
    \leq c_1\eps^{1/2}\|f\|.
\end{eqnarray}

Next, since $W_2^2(\Real\setminus \{0\})\subset C^1(\Real\setminus \{0\})$ by the Sobolev embedding theorem, one gets a constant $c_2>0$ such that
\(    
    \|g\|_{C^1(\Real\setminus \{0\})}\leq c_2\|g\|_{W_2^2(\Real\setminus \{0\})}
\)    
for all $g \in W_2^2(\Real\setminus \{0\})$.
\rd{Using the Sobolev embedding theorem on the interval $[-1,1]$ and properties of solutions to the problem~\eref{WepsProblem}, one arrives at the estimates}
\begin{eqnarray}\label{EstForWeps}
\rd{  \|w_\eps\|_{C^1[-1,1]}\leq c_2 \|w_\eps\|_{W_2^2(-1,1)}\leq c_3\|f(\eps\cdot)\|_{L_2(-1,1)}\leq c_3\eps^{-1/2} \|f\|}
\end{eqnarray}
\rd{with some constant $c_3$ independent of~$f$.}
Therefore, taking into account the equality $y(0+) = \theta y(0-)$, we see that the jump of $z_\eps$ at~$x=\eps$ can be bounded as
\begin{eqnarray}
 \fl \bigr|[z_\eps]_{\eps}\bigl|
        =\bigr|y(\eps)-\theta y(0-)-\eps y'(0-)v(1)\bl{-\eps^2 w_\eps(1)}\bigl|\nonumber\\
        \leq \bigr|y(\eps)-y(0+)\bigl| +\eps \bigl|y'(0-)\bigr|\bigl|v(1)\bigr|\bl{+\eps^2\bigl| w_\eps(1)\bigr|} \nonumber\\
        \leq (\eps^{1/2}+c_2|v(1)|\,\eps) \|y\|_{W_2^2(\Real\setminus \{0\})}\bl{+c_3\eps^{3/2}\|f\|}
        \leq \bl{c_4}\eps^{1/2}\|f\| \label{JumpZAtEps}
\end{eqnarray}
for some \bl{$c_4$} independent of~$f$ and $\eps\in(0,1)$.
Similarly, keeping in mind that $v'(1)= \theta^{-1}$ and  $y'(0+)=\theta^{-1}y'(0-)$ \rd{and using~\eref{SinvEst} and \eref{EstForWeps}}, we get
\begin{eqnarray}
 \fl  \bigr|[z'_\eps]_{\eps}\bigl|
    =\bigr|y'(\eps)-y'(0-)v'(1)\bl{-\eps w'_\eps(1)}\bigr|
    \rd{\leq}\bigl|y'(\eps)-y'(0+)\bigr|\bl{+\eps \bigl|w'_\eps(1)\bigr|}\nonumber\\
    \leq \eps^{1/2}\|y\|_{W_2^2(\Real\setminus \{0\})}\rd{+c_3\eps^{1/2}\|f\|
    \leq c_5\eps^{1/2}\|f\|}.\label{JumpZ'AtEps}
\end{eqnarray}

\begin{figure}[t]
\begin{center}
\includegraphics[scale=1]{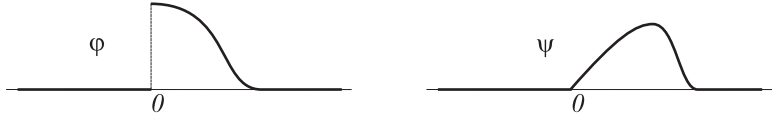}
\end{center}\caption{\label{Jump} Functions with the prescribed jumps at the origin}
\end{figure}

Let us introduce functions $\varphi$ and $\psi$ as on \fref{Jump}  that are smooth outside the origin, have compact supports contained in $[0,\infty)$, and have the prescribed jumps
$[\varphi]_0=1$, $[\varphi']_0=0$ and $[\psi]_0=0$, $[\psi']_0=1$.
Set
\begin{equation}\label{Zeta}
\fl\zeta_\eps(x)=
    [z_\eps]_{-\eps}\,\varphi(-x-\eps)-[z'_\eps]_{-\eps}\,\psi(-x-\eps)
        -[z_\eps]_{\eps}\,\varphi(x-\eps)-[z'_\eps]_{\eps}\,\psi(x-\eps);
\end{equation}
then $\zeta_\eps =0$ on~$(-\eps,\eps)$ and, in view of~\rd{\eref{JumpZAt-Eps}--\eref{JumpZ'At-Eps} and \eref{JumpZAtEps}--\eref{JumpZ'AtEps}},
\begin{equation}\label{ZetaEstim}
    \max_{|x|>\eps}\left|\zeta^{(k)}_\eps(x)\right|\leq \bl{c_6}\eps^{1/2}\|f\|
\end{equation}
for some $\bl{c_6}>0$, $k=0,1,2$, and all~$\eps>0$. Clearly, the function $y_\eps:= z_\eps+\zeta_\eps$ is continuous on~$\Real$ along with its derivative and belongs to $W_2^2(\Real) = \dmn S_\eps$. Observe that $y_\eps=y$ for $|x|$ large enough; more exactly, we have
\begin{equation}\label{yeps}
    y_\eps =  y +q_\eps=(S(\theta)-k^2)^{-1}f+q_\eps
\end{equation}
with
\[
  \fl  q_\eps(x):= -\chi(x/\eps)y(x)+y(0-)u(x/\eps)+\eps y'(0-) v(x/\eps)\bl{+\eps^2 w_\eps(x/\eps)}+\zeta_\eps(x)
\]
of compact support. The first \rd{four} summands above are functions bounded uniformly in~$\eps$ and have support contained in~$[-\eps,\eps]$, and the \rd{last} one is small due to~\eref{ZetaEstim}. It follows from the above estimates that, for a suitable $c_5$ independent of~$f$ and $\eps\in(0,1)$,
\begin{equation}\label{QepsEst}
    \|q_\eps\|\leq \bl{c_7}\eps^{1/2}\|f\|.
\end{equation}

We are now in a position to prove \eref{SepsEst}.
If $\abs{x}>\eps$, then
\begin{equation}\label{ActionB}
     (S_\eps-k^2)y_\eps= \Bigl(-\frac{\rmd^2}{\rmd x^2} -k^2\Bigr)(y + q_\eps)= f-\zeta''_\eps-k^2\zeta_\eps.
\end{equation}
If $|x|<\eps$, then $\zeta_\eps(x)=0$, and so $y_\eps(x)=z_\eps(x)= y(0-)u(\tfrac{x}{\eps})+\eps y'(0-) v(\tfrac{x}{\eps})\rd{+\eps^2 w_\eps(\tfrac{x}{\eps})}$. \rd{Thus,}
\begin{eqnarray}
\fl (S_\eps-k^2)y_\eps
    =\eps^{-2}y(0-)\left\{-u''(\tfrac{x}{\eps})
            +V(\tfrac{x}{\eps})u(\tfrac{x}{\eps})\right\}
    +\eps^{-1}y'(0-) \left\{-v''(\tfrac{x}{\eps})
        + V(\tfrac{x}{\eps})v(\tfrac{x}{\eps})\right\}\nonumber\\
   \bl{+ \left\{-w_\eps''(\tfrac{x}{\eps})
        + V(\tfrac{x}{\eps})w_\eps(\tfrac{x}{\eps})\right\} }
    - k^2 \rd{\chi(\tfrac{x}{\eps})y_\eps(x)}       \label{ActionS} \\
    =   \rd{\chi(\tfrac{x}{\eps})[f(x)- k^2 y_\eps(x)]}, \nonumber
\end{eqnarray}
since both $u$, $v$ are solutions to equation \eref{EqV} \bl{and $w_\eps$ is a solution to \eref{WepsProblem}.}  Therefore,
\begin{equation}\label{ActionG}
    (S_\eps-k^2)y_\eps=f+r_\eps,
\end{equation}
where
\begin{eqnarray*}
   r_\eps(x)&= -\zeta''_\eps(x)
                    - \rd{k^2\left[\zeta_\eps(x)+\chi(x/\eps) y_\eps(x/\eps)\right]}\\
        &= -\zeta''_\eps(x) -\rd{k^2\left[q_\eps(x)+\chi(x/\eps)y(x)\right].}
\end{eqnarray*}
Relations~\eref{SinvEst}, \eref{ZetaEstim}, and \eref{QepsEst} now yield the estimate
\begin{equation}\label{RepsESt}
\|r_\eps\|\leq c \eps^{1/2}\|f\|
\end{equation}
for a suitable $c\ge \bl{c_7}$, and \eref{SepsEst} is proved.\hfill $\square$
\medskip

\begin{thm}\label{thm:reson}
  Assume that $\lambda=0$ is an eigenvalue of the operator $\mathcal{N}$ with an eigenfunction $u$ and set $\theta:=u(1)/u(-1)$. Then $S_\eps$ converge to
  $S(\theta)$ as $\eps\to 0$ in the norm resolvent sense.
\end{thm}
\smallskip

\noindent
\textbf{Proof.} Fix $k^2 \in \mathbb{C}\setminus\mathbb{R}$. We conclude from \eref{ActionG} and \eref{yeps} that
\begin{equation*}
(S_\eps-k^2)^{-1}f=y_\eps-(S_\eps-k^2)^{-1}r_\eps=(S(\theta)-k^2)^{-1}f+q_\eps-(S_\eps-k^2)^{-1}r_\eps
\end{equation*}
for each $f\in L_2(\Real)$. This gives   $(S_\eps-k^2)^{-1}f-(S(\theta)-k^2)^{-1}f=q_\eps-(S_\eps-k^2)^{-1}r_\eps$,
so that, by \eref{QepsEst} and \eref{RepsESt},
  \begin{eqnarray*}
  \fl  \|(S_\eps-k^2)^{-1}f-(S(\theta)-k^2)^{-1}f\| \leq  \|q_\eps\|+\|(S_\eps-k^2)^{-1}r_\eps\|\\
  \qquad\qquad\qquad\quad\qquad\qquad\qquad \leq
    \|q_\eps\|+|{\Im}\,k^2|^{-1}\|r_\eps\|\leq \bl{c_8}\eps^{1/2}\|f\|
  \end{eqnarray*}
for a suitable \bl{$c_8$} independent of $f$ and $\eps$. The theorem is proved. \hfill $\square$

\section{Non-Resonant Case}\label{sec:non-res}
Now we study the non-resonant case when $\lambda=0$ is not an eigenvalue of $\mathcal{N}$.
Recall that $S_-\oplus S_+$ denotes the direct sum  of the unperturbed half-line Schr\"odinger operators $S_\pm= -\rmd^2/\rmd x^2$ on~$\Real_\pm$ subject to the Dirichlet boundary condition at $x=0$;
we shall prove that $S_-\oplus S_+$ is the limit of $S_\eps$ as $\eps \to0$.

The above proof still works for the non-resonant case after we have slightly changed the corrector $q_\eps$ of Lemma~\ref{LemmaResonant}. Let $y=(S_-\oplus S_+-k^2)^{-1}f$ for some $f\in L_2(\Real)$ and let $w$ be a unique solution of the problem
  $$ -w''+V(\xi)w=0,\quad \xi\in(-1,1), \qquad w'(-1)=y'(0-),\quad w'(1)=y'(0+).$$
It is well known that $w$ obeys the a priori estimate $ \|w\|_{W_2^2}\leq c_1 \bigr(|y'(0-)|+|y'(0+)|\bigl)$
for some constant $c_1$ independent of~$y'(0\pm)$. Combining it with the reasoning of the previous section, we find positive constants $c_2$, $c_3$, and $c_4$ such that the following inequalities hold:
\begin{equation}\label{VestNR}
   \|w\|_{C([-1,1])} \leq c_2 \|w\|_{W_2^2} \leq c_3 \|y\|_{W_2^2(\Real\setminus \{0\})}\leq c_4\|f\|.
\end{equation}

We now extend the function~$w$ to the whole line by zero and introduce the function
\begin{equation*}
z_\eps(x)=(1-\chi(x/\eps))y(x)+\eps w(x/\eps)\bl{+\eps^2 w_\eps(x/\eps),}
\end{equation*}
\bl{where $w_\eps$ is a solution to \eref{WepsProblem} as above.}
The jumps of $z_\eps$ and $z_\eps'$ at the points $x=\pm\eps$  converge to 0 as $\eps\to 0$ uniformly with respect to the $L_2$-norm of $f$. In fact, taking into account the conditions $y(0-)=0= y(0+)$ we obtain
\begin{eqnarray*}
\fl  &\bigr|[z_\eps]_{-\eps}\bigl|=\bigr|\eps w(-1)-y(-\eps)\bigl|\leq \eps|w(-1)|+\int^0_{-\eps}|y'(t)|\,\rmd t\leq c_5\eps^{1/2}\|f\|,\\
\fl    &\bigr|[z_\eps]_{\eps}\bigl|=\bigr|y(\eps)-\eps w(1)\bl{-\eps^2w_\eps(1)}\bigl|\leq \int_0^{\eps}|y'(t)|\,\rmd t+\eps|w(1)|\bl{+\eps^2|w_\eps(1)|}\leq c_6\eps^{1/2}\|f\|,
\end{eqnarray*}
due to \eref{VestNR}.
Next, the jumps $[z'_\eps]_{\pm\eps}$  can be estimated  as in \eref{JumpZ'At-Eps}, \eref{JumpZ'AtEps} above.

Let us introduce the function $y_\eps=z_\eps+\zeta_\eps$, where $\zeta_\eps$ is defined as in \eref{Zeta}. This function belongs to $\dmn S_\eps$ and can be written in the form
\begin{equation*}
y_\eps=(S_-\oplus S_+-k^2)^{-1}f+q_\eps
\end{equation*}
with  $q_\eps(x)=-\chi(x/\eps)y(x)+\eps w(x/\eps)\bl{+\eps^2 w_\eps(x/\eps)}+\zeta_\eps(x)$.
The $L_2$-norm of the function $q_\eps$  can be estimated by $c\eps^{1/2}\|f\|$.
By calculations similar to those in \eref{ActionB} and \eref{ActionS}, we establish that
\begin{equation*}
(S_\eps-k^2)y_\eps=f+r_\eps,
\end{equation*}
where $r_\eps(x)=   -\zeta''_\eps(x) -\rd{k^2\left[q_\eps(x)+\chi(x/\eps)y(x)\right]}$.
As above, $\|r_\eps\|\leq c \eps^{1/2}\|f\|$ for some constant~$c$ independent of~$f$.

We now get the following theorem, whose proof is analogous to that of Theorem~\ref{thm:reson}.

\begin{thm}\label{thm:nonreson}
  If $\lambda=0$ is not an eigenvalue of the operator $\mathcal{N}$, then the family~$S_\eps$ of~\eref{eq:intr.Seps} converges in the norm resolvent sense as $\eps\to 0$ to~$S_-\oplus S_+$.
\end{thm}

\section{What is the right Hamiltonian with a $\delta'$-potential?}\label{sec:delta'}

Let us consider the formal Schr\"{o}dinger operators
\begin{equation*}
H_\alpha = -\frac{\rmd^2}{\rmd x^2} + \alpha \delta'(x),
\end{equation*}
whose potentials contain the derivative of the Dirac delta function. Here $\alpha$ is a strength interaction parameter or a coupling constant taking values in~$\Real$.
Equation $-v''+\alpha \delta'(x)v=\lambda v$ has no solutions in the space of distributions, except for the trivial one.
Nevertheless, the formal Hamiltonian $H_\alpha$ may be defined
in terms of the equation $-v''=\lambda v$, for $x\neq0$, and appropriate
boundary conditions at the origin.
The main question is therefore how to choose these boundary conditions or, in other words, how to choose a proper self-adjoint extension of the so-called minimal operator
$$
  S'_{0}=-\frac{\rmd^2}{\rmd x^2}, \quad \dmn S_0'=\{ g \in W_2^2(\Real) \mid g(0) =g'(0)=0\}.
$$
The minimal operator~$S_0'$ is symmetric and has deficiency indices~$(2,2)$; therefore, its self-adjoint extensions form a four-parametric family and there are different possibilities to define $H_\alpha$.

\begin{figure}[t]
\begin{center}
\includegraphics[scale=0.2]{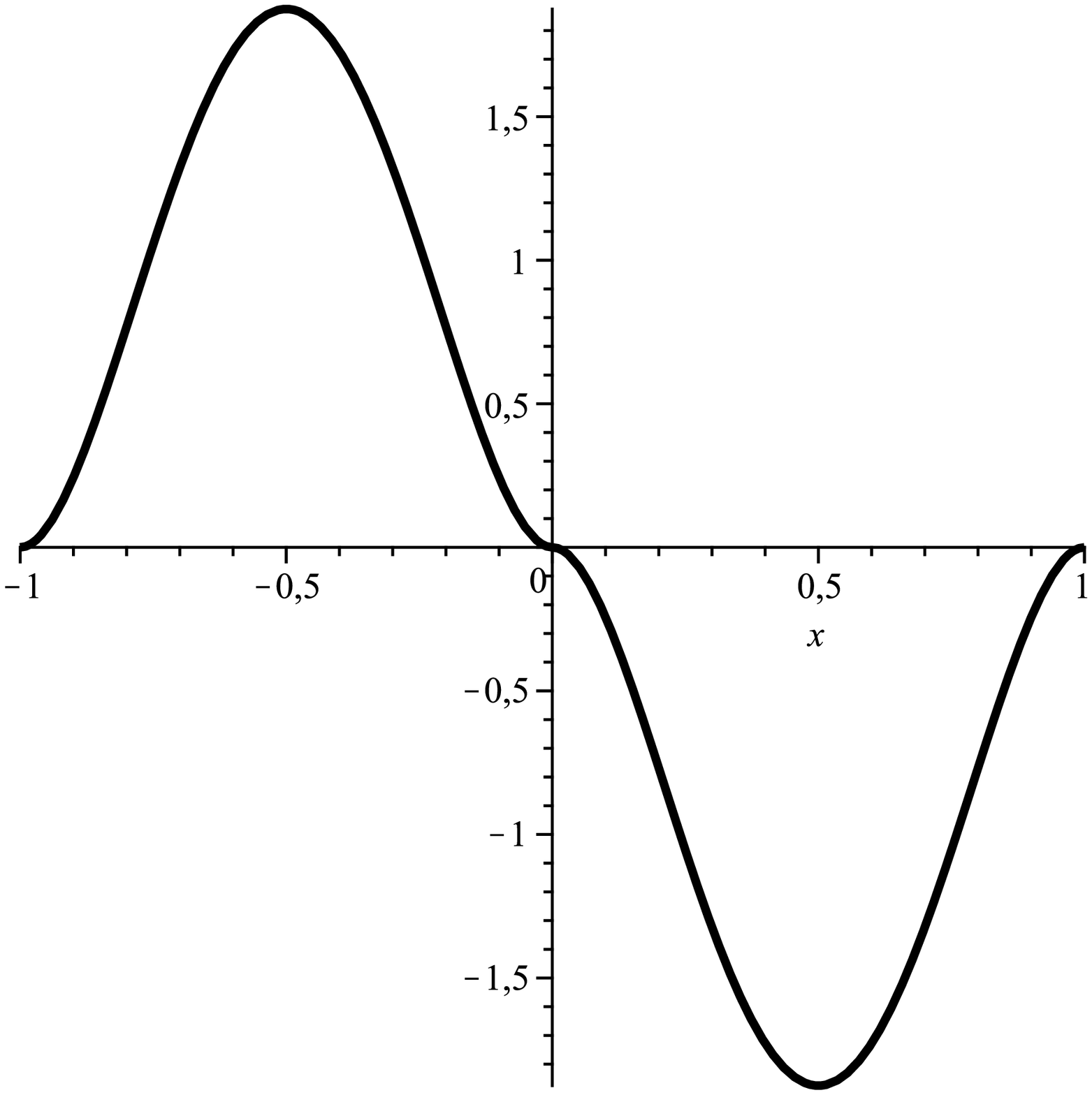}\quad\qquad\includegraphics[scale=0.2]{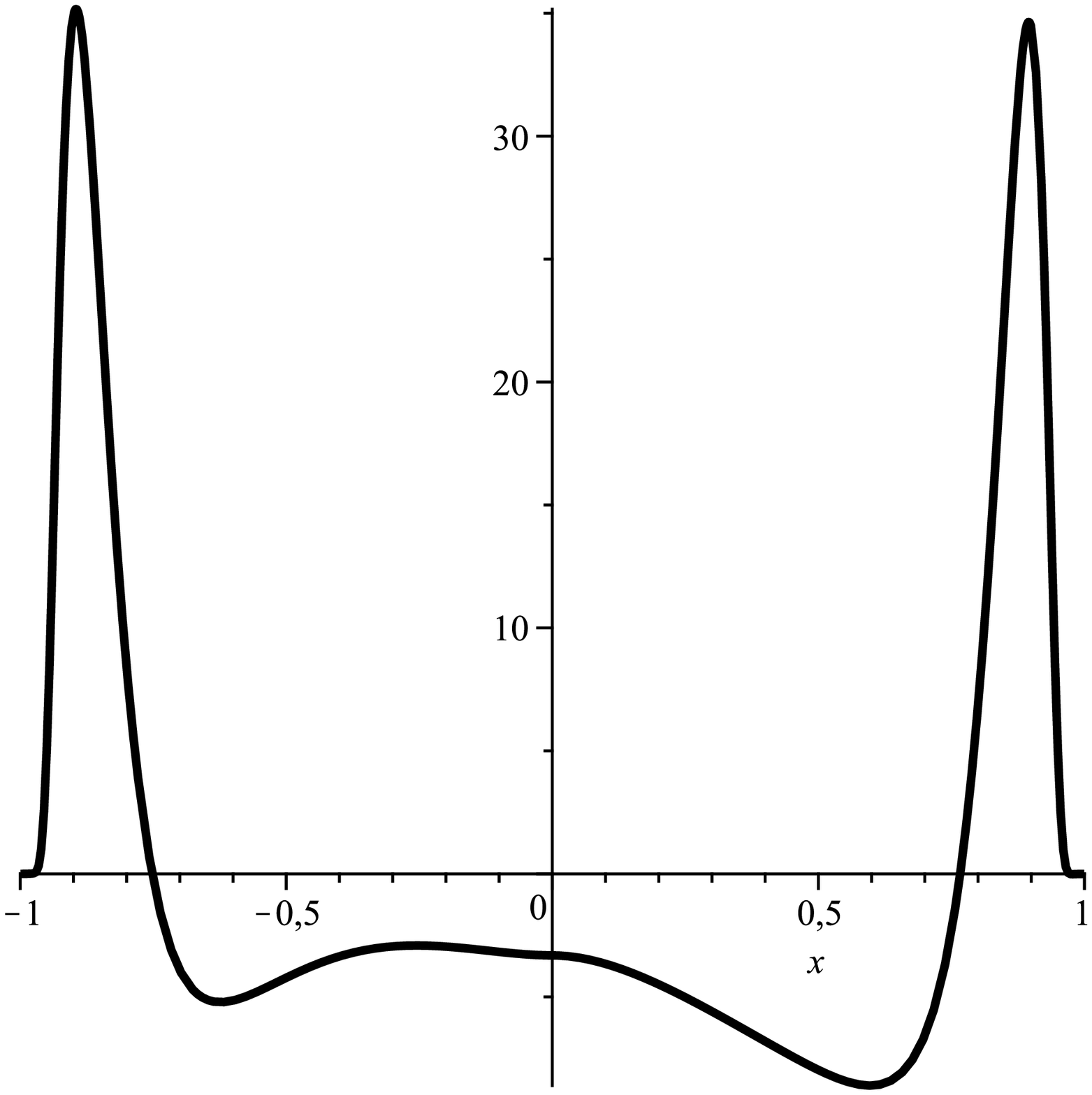}\quad\qquad\includegraphics[scale=0.2]{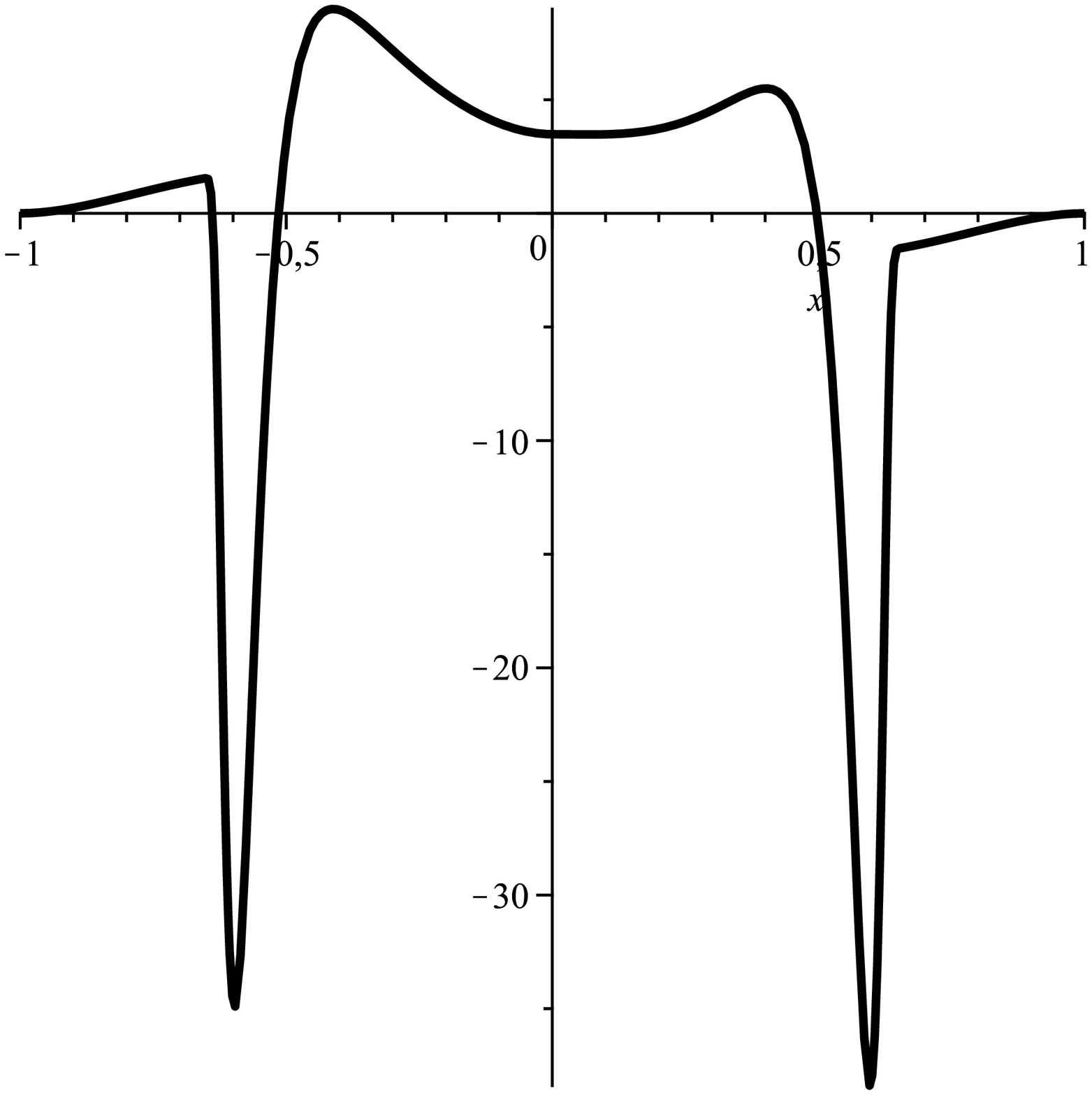}
\end{center}\caption{Odd and more exotic $\delta'$-like potentials}\label{figDP}
\end{figure}

An alternative way is to realize the operator $H_\alpha$ as the limit of Hamiltonians with regularized potentials
\begin{equation}\label{Heps}
   H_\eps(\alpha,\Psi)= -\frac{\rmd^2}{\rmd x^2}+\frac{\alpha}{\eps^2}\Psi(\eps^{-1}x)
\end{equation}
with $\eps$ being a regularization parameter. Suppose that $\Psi\in C_0^\infty(\Real)$ and $\supp \Psi=[-1,1]$.
It is easy to check that $\eps^{-2}\Psi(\eps^{-1}x)\to \delta'(x)$ in the sense of distributions as $\eps\to0$ iff
\begin{equation}\label{Moments}
    \int_\Real \Psi(\xi)\,\rmd\xi=0\quad\mbox{and}\quad\int_\Real \xi \Psi(\xi)\,\rmd\xi=-1;
\end{equation}
in this case, we call $\Psi$ a $\delta'$\emph{-like potential}. Plots of some $\delta'$-like potentials are shown on \fref{figDP}.

Let us introduce two spectral characteristics of a $\delta'$-like potential~$\Psi$. We denote by~$\Sigma_\Psi$
the set of $\alpha$-eigenvalues of the problem
\begin{equation}\label{ResonantProblem}
    -w''+\alpha\Psi(\xi)\,w=0,\quad \xi\in(-1,1),\qquad w'(-1)=w'(1)=0
\end{equation}
and call $\Sigma_\Psi$ the \emph{resonant set} of $\Psi$.
It is well known \cite{CurgusLangerJDE1989} that the $\alpha$-spectrum of the problem~\eref{ResonantProblem} is discrete and real and  consists of simple non-zero eigenvalues and the  geometrically simple eigenvalue $\alpha=0$; to the latter, there correspond a constant eigenfunction and an adjoint function. Also, $\Sigma_\Psi$ possesses two accumulation points at~$\pm\infty$ since the function $\Psi$ changes sign.
It is obvious that $\alpha$ belongs to the resonant set if the potential $\alpha\Psi$ is resonant in the sense of the definition given in the introduction.
For any non-zero $\alpha\in \Sigma_\Psi$, we put
\begin{equation*}
    \theta_\Psi(\alpha)=\frac{w_\alpha(1)}{w_\alpha(-1)},
\end{equation*}
where $w_\alpha$ is an eigenfunction corresponding to $\alpha$. The ratio is properly defined since the value $w_\alpha(-1)$ is different from $0$. Besides, $\theta_\Psi(\alpha)$ is always real and does not depend on the choice of the eigen\-function. We call $\theta_\Psi\colon \Sigma_\Psi\to \Real$ \textit{the coupling function} of $\Psi$.

Let us introduce the one parameter family of self-adjoint extensions of $S_0'$:
\begin{equation*}
    H(\alpha,\Psi)=
   \left\{
     \begin{array}{ll}
       S_-\oplus S_+, & \hbox{if $\alpha\not\in \Sigma_\Psi$,} \\
       S(\theta_\Psi(\alpha)), & \hbox{if $\alpha\in\Sigma_\Psi$.}
     \end{array}
   \right.
\end{equation*}
Hence $H(\alpha,\Psi)$ is a \emph{connected} self-adjoint extension with the domain given by~\eref{DomS(theta)} with $\theta=\theta_\Psi(\alpha)$
for the resonant coupling constants $\alpha$; otherwise, $H(\alpha,\Psi)$ is a \emph{separated} extension, namely the direct sum of unperturbed half-line Schr\"odinger operators subject to the Dirichlet boundary conditions.

\begin{thm}
  Assume $\Psi$ is a $\delta'$-like potential. Then the family of Hamiltonians $H_\eps(\alpha,\Psi)$
 converges in the norm resolvent sense as $\eps\to 0$ to $H(\alpha,\Psi)$.
\end{thm}
\smallskip

\noindent
\textbf{Proof.} This theorem is  a reformulation of  theorems~\ref{thm:reson} and \ref{thm:nonreson} in terms of the po\-ten\-ti\-al~$\alpha \Psi$. Clearly, $\lambda=0$ is an eigenvalue of the operator $\mathcal{N}$ with potential $V=\alpha \Psi$ if and only if the coupling constant $\alpha$ belongs to the resonant set $\Sigma_\Psi$.
 \hfill $\square$
\medskip

Note that $C^\infty$-smoothness of $\Psi$ is not essential and the theorem also remains valid for $\Psi\in L_1(\Real)$. Moreover, the theorem holds not only for $\delta'$-like potentials~$\Psi$, but also for arbitrary real-valued $\Psi$ of compact support not necessarily satisfying~\eref{Moments}.

Provided that $\Psi$ is the ``shape'' of a $\delta'$-like short range potential in the physical model under consideration, we can define the limiting Schr\"odinger operator with a po\-ten\-ti\-al~$\alpha \delta'$ as the operator $H(\alpha,\Psi)$. As a result, we find an interesting feature of the $\delta'$-potentials: the appropriate solvable model is not unique and crucially depends on the way in which the derivative of the Dirac delta-function is approximated in the weak topology. The shape $\Psi$ is a ``hidden'' parameter in the conventional formulation of the problem on the right definition of the Hamiltonians with $\delta'$-like potentials.


\section{Resonant phenomenon  in  transmission through a $\delta'$-potential}\label{sec:transm}

It is important to emphasize that the scattering properties of the point dipole also depend on the way in which the zero-range limit is realized. We shall show that the transmission coefficient depends on
the intensity $\alpha$ and the regularizing profile~$\Psi$ in such a way that for all values of $\alpha$, the $\delta'$-barrier is completely impenetrable except for the set~$\Sigma_\Psi$ of resonant values, at which there occurs a partial transmission fading away as $|\alpha|$ becomes larger.

First, let us discuss stationary scattering associated with the Hamiltonians
$H(\alpha,\Psi)$ and $-\frac{d^2}{dx^2}$.
We need to consider only the case when the interaction parameter $\alpha$ belongs to the resonant set $\Sigma_\Psi$. Consider the incoming monochromatic wave $e^{ikx}$ with $k>0$ coming from the
left. Then the corresponding wave function has the form
\begin{equation*}
    \psi(x,k)=\left\{
                \begin{array}{ll}
                  e^{ikx} + R\,e^{-ikx} & \hbox{if $x<0$,} \\
                  \phantom{e^{ikx}+\,}\  T\,e^{ikx} & \hbox{if $x>0$.}
                \end{array}
              \right.
\end{equation*}
Here $R$ and $T$ are respectively the reflection and transmission coefficients. As far as  $H(\alpha,\Psi)=S(\theta_\Psi(\alpha))$, the matching
conditions at the origin
\begin{equation*}
\left(
\begin{array}{c}
  \psi(0+,k) \\
  \psi'(0+,k)
\end{array}
\right)=
\left(
\begin{array}{cc}
 \theta_\Psi(\alpha) & 0\\
  0 & \theta_\Psi(\alpha)^{-1}
\end{array}
\right)
\left(
\begin{array}{c}
  \psi(0-,k) \\
  \psi'(0-,k)
\end{array}
\right)
\end{equation*}
clearly yield
\begin{equation*}
\left(
\begin{array}{c}
  T \\
  ikT
\end{array}
\right)=
\left(
\begin{array}{cc}
 \theta_\Psi(\alpha) & 0\\
  0 & \theta_\Psi(\alpha)^{-1}
\end{array}
\right)
\left(
\begin{array}{c}
  1+R \\
  ik(1-R)
\end{array}
\right).
\end{equation*}
Then one obtains the values of the reflection and transmission coefficients that can be expressed via the coupling function as follows (cf.\ also\cite{AlbeverioCacciapuotiFinco:2007}):
\begin{equation}\label{LimitCoeffs}
\fl    R_\Psi(\alpha)=\left\{
                     \begin{array}{ll}
                     \displaystyle  \frac{1-\theta_\Psi^2(\alpha)}{1+\theta_\Psi^2(\alpha)} & \hbox{if $\alpha\in\Sigma_\Psi$,} \\
                       \qquad -1 & \hbox{otherwise,}
                     \end{array}
                   \right.\qquad
T_\Psi(\alpha)=\left\{
                 \begin{array}{ll}
                  \displaystyle   \frac{2\theta_\Psi(\alpha)}{1+\theta_\Psi^2(\alpha)} & \hbox{if $\alpha\in\Sigma_\Psi$,} \\
                   \qquad 0 & \hbox{otherwise.}
                 \end{array}
               \right.
\end{equation}
We emphasize that the reflection and transmission coefficients do not depend on $k$.

Next we investigate stationary scattering for the Hamiltonians $H_\eps(\alpha,\Psi)$ and $-\frac{d^2}{dx^2}$ and prove that the scattering amplitude converges as $\eps\to0$ to that for the limiting Hamiltonian $H(\alpha,\Psi)$. The analysis below basically follows a more general approach of~\cite{Manko:2009}.

We are looking for the positive-energy solution of the equation with a $\delta'$-like potential
\begin{equation*}
-y''+\alpha\eps^{-2}\Psi(\eps^{-1}x)y=k^2y, \quad x\in\Real,
\end{equation*}
 given in the form
\begin{equation*}
\psi_\eps(x,k,\alpha)=\left\{
                     \begin{array}{ll}
                       e^{ikx}+R\,e^{-ikx} & \hbox{if $x<-\eps$,} \\
   A\,u(\eps^{-1}x,\eps k,\alpha)+B\,v(\eps^{-1}x,\eps k,\alpha) & \hbox{if $|x|<\eps$,} \\
                       \phantom{e^{ikx}+\,}T\,e^{ikx} & \hbox{if $x>\eps$.}
                     \end{array}
                   \right.
\end{equation*}
Here $u=u(\xi,\varkappa,\alpha)$ and $v=v(\xi,\varkappa,\alpha)$ are the solutions of the equation
\begin{equation}\label{AuxlProblem}
-w''+\alpha\Psi(\xi)w=\varkappa^2 w, \quad \xi\in(-1,1)
\end{equation}
subject to the  initial conditions $u(-1,\varkappa,\alpha)=1$, $u'(-1,\varkappa,\alpha)=0$ and $v(-1,\varkappa,\alpha)=0$, $v'(-1,\varkappa,\alpha)=1$ respectively.
The unknown coefficients $R$, $A$, $B$, and $T$ can be found from the linear system (set~$\varkappa:=\eps k$)
\begin{equation*}
     \left( \begin{array}{cccc}
-e^{i\varkappa}& 1& 0& 0 \\
     i\varkappa e^{i\varkappa}& 0& 1&0\\
     0&u(1,\varkappa,\alpha)&v(1,\varkappa,\alpha)&-e^{i\varkappa}\\
     0&u'(1,\varkappa,\alpha)&v'(1,\varkappa,\alpha)&-i\varkappa e^{i\varkappa}
       \end{array}\right)
     \left(\begin{array}{c} R\\A\\B\\T\end{array}\right)=
    \left(\begin{array}{c} e^{-i\varkappa} \\i\varkappa e^{-i\varkappa}\\0\\0\end{array}\right)
\end{equation*}
obtained by matching the solution and its first derivative at the points $x=\pm\eps$.
The system determinant admits the asymptotic expansion
\begin{equation}\label{Asymptotyc}
   \Delta(\varkappa,\alpha)=u'(1,0,\alpha)+i\varkappa\, q(\alpha)+O(\varkappa^2),\quad\varkappa\to 0,
\end{equation}
where $q(\alpha)=2u'(1;0,\alpha)-u(1;0,\alpha)- v'(1;0,\alpha)$. By Cramer's rule, one obtains
\begin{eqnarray}\label{Reps}
R_\Psi(\varkappa,\alpha)=
    \frac{-u'(1;0,\alpha)
    +i\varkappa[u(1;0,\alpha)-v'(1;0,\alpha)]}{u'(1;0,\alpha)
    +i\varkappa\,q(\alpha)}+O(\varkappa^2),\\ \label{Teps}
T_\Psi(\varkappa,\alpha)=
    \frac{-2i\varkappa}{u'(1;0,\alpha)+i\varkappa\,q(\alpha)}+O(\varkappa^2)
\end{eqnarray}
as $\varkappa\to0$.
Here we use the identity $u(1,\varkappa,\alpha)v'(1,\varkappa,\alpha)-u'(1,\varkappa,\alpha)v(1,\varkappa,\alpha)=1$ that follows from the constancy in $\xi$ of the Wronskian of $u$ and $v$.

\begin{thm}
  For each $k>0$  and  $\alpha\in\Real$ the reflection and transmission coefficients
$R_\Psi(\eps k,\alpha)$ and\, $T_\Psi(\eps k,\alpha)$ converge towards $R_\Psi(\alpha)$
and\, $T_\Psi(\alpha)$
as $\eps\to 0$ respectively, where the limit values are given by \eref{LimitCoeffs}.
\end{thm}
\smallskip

\noindent
\textbf{Proof.} \textit{The non-resonant case.} Since the equation~\eref{AuxlProblem} for $\varkappa=0$ coincides with~\eref{ResonantProblem} and $\alpha$ is not a resonant coupling constant, we conclude that $u'(1;0,\alpha)$ is different from $0$. From \eref{Reps} and \eref{Teps}, it immediately follows that $R_\Psi(\eps k,\alpha)= -1+O(\eps k)$ and $T_\Psi(\eps k,\alpha)=O(\eps k)$ as $\eps\to 0$.

\textit{The resonant case.} If $\alpha\in \Sigma_\Psi$, then $u'(1;0,\alpha)=0$. Hence $u$ is an eigenfunction of \eref{ResonantProblem} with eigenvalue $\alpha$. Next, $u(1;0,\alpha)=\theta_\Psi(\alpha)$ by the definition of the coupling function and $v'(1;0,\alpha)=\theta_\Psi(\alpha)^{-1}$ by the Lagrange identity, which yields  the relation $q(\alpha)=-(\theta_\Psi(\alpha)^{-1}+\theta_\Psi(\alpha))$. Therefore
\begin{eqnarray*}
     R_\Psi(\eps k,\alpha)=
    \frac{\theta_\Psi(\alpha)^{-1}-\theta_\Psi(\alpha)}{\theta_\Psi(\alpha)^{-1}+\theta_\Psi(\alpha)}+O(\eps^2 k^2)=
\frac{1-\theta_\Psi^2(\alpha)}{1+\theta_\Psi^2(\alpha)}+O(\eps^2 k^2),\\
T_\Psi(\eps k,\alpha)=
    \frac{2}{\theta_\Psi(\alpha)^{-1}+\theta_\Psi(\alpha)}+O(\eps^2 k^2)=\frac{2\theta_\Psi(\alpha)}{1+\theta_\Psi^2(\alpha)}+O(\eps^2 k^2)
\end{eqnarray*}
as $\eps\to 0$, and the proof is complete.
 \hfill $\square$
\smallskip


\section{An example}\label{sec:example}

Let us assume that the shape of a short-range potential in an actual model can be approximately described as
\begin{equation*}
    \Psi(x)=\left\{
                     \begin{array}{ll}
                       -6x(x+1) & \hbox{ if $x\in [-1,0]$,} \\
   \phantom{-}6x(x-1) & \hbox{ if $x\in [0,1]$,} \\
                      \qquad 0 & \hbox{ otherwise.}
                     \end{array}
                   \right.
\end{equation*}
The function~$\Psi$ is a $\delta'$-like potential, i.e.\ it satisfies conditions \eref{Moments}.
As shown above, the best choice of the solvable model corresponding to the family of Hamiltonians
\[
    -\frac{d^2}{dx^2}+\frac{\alpha}{\eps^2}\Psi(\eps^{-1}x)
\]
in the zero-range limit is given by the operator~$H(\alpha,\Psi)$. Certainly, the resonant set~$\Sigma_\Psi$ and the coupling function $\theta_\Psi$
to be found are specific to the given shape~$\Psi$.

\begin{table}[t]\label{table}
\caption{\label{tabone}Resonant intensities, coupling function and transmission probabilities}
\begin{indented}
\lineup
\item[]\begin{tabular}{*{4}{ccccc}}
\br
$\alpha$&&$\theta_\Psi(\alpha)$&&$|T_\Psi(\alpha)|^2$\cr
\mr
0&& 1 &&$1$\cr
\0 18.1747&& \0\0 -54.9385&&0.00132\cr
\0 57.1490&&\0 1352.8032&&  $0.219\cdot 10^{-5}$\cr
117.4863&&-32156.4597&&     $0.387\cdot 10^{-8}$ \cr
199.1756&&755821.4703&&     $0.704\cdot 10^{-11}$\cr
\br
\end{tabular}
\end{indented}
\end{table}

Since the function $\Psi$ is odd, the set $\Sigma_\Psi\subset \Real$ is symmetric with respect to the origin, namely if $\alpha$ is an eigenvalue of \eref{ResonantProblem} with eigenfunction $w_\alpha$, then $-\alpha$ is also an eigenvalue of \eref{ResonantProblem} with eigenfunction $w_{-\alpha}(\xi)=w_{\alpha}(-\xi)$.
Furthermore we  conclude from this that
\begin{equation*}
    \theta_\Psi(-\alpha)= \frac{w_{-\alpha}(1)}{w_{-\alpha}(-1)}=\frac{w_\alpha(-1)}{w_\alpha(1)} =\frac{1}{  \theta_\Psi(\alpha)};
\end{equation*}
hence, as follows from~\eref{LimitCoeffs}, $|R_\Psi(-\alpha)|^2=|R_\Psi(\alpha)|^2$ and $|T_\Psi(-\alpha)|^2=|T_\Psi(\alpha)|^2$.

Table~\ref{tabone} lists the first five nonnegative resonant values of~$\alpha$ (numerically computed using \textit{Maple}) and the corresponding values of the coupling $\theta_\Psi(\alpha)$ and the squared transmission coefficient~$|T_\Psi(\alpha)|^2$. We note that the latter decays very fast and conjecture that this will be observed for all $\delta'$-like profiles~$\Psi$.

\ack The authors are grateful to S~Albeverio and C~Cacciapuoti for bringing to their attention the papers~\cite{AlbeverioCacciapuotiFinco:2007,CacciapuotiExner:2007,Seba:1985} and stimulating discussions and to S Man'ko for providing them with the asymptotic analysis of the transmission and reflection coefficients of~\cite{Manko:2009} in section~\ref{sec:transm}. They also thank the anonymous referees for careful reading of the manuscript and valuable remarks and suggestions. The research of RH was partially supported by Deutsche Forschungsgemeinschaft under project 436 UKR 113/84.


\section*{References}
\begin{thebibliography}{10}

\bibitem{AlbeverioCacciapuotiFinco:2007}
    Albeverio S, Cacciapuoti C, and Finco D 2007
    \emph{J. Math. Phys.}  \textbf{48} no.~3, 032103, 21 pp

\bibitem{Albeverio2edition}
Albeverio S,  Gesztesy F,  H{\o}egh-Krohn R and  Holden H 2005
{\it Solvable Models in Quantum Mechanics. With an Appendix by Pavel Exner. 2nd revised edn}
(Providence, RI: AMS Chelsea Publishing) p~488

\bibitem{AlbeverioKoshmanenko:1999a}
    Albeverio S and Koshmanenko V 1999
    {\it Potential Anal.} {\bf11} no.~3, 279

\bibitem{AlbeverioKoshmanenko}
    Albeverio S and Koshmanenko V 1999   \emph{J. Funct. Anal.}  \textbf{169} no.~1, 32

\bibitem{AlbeverioKurasov}
    Albeverio S and Kurasov P 1999
    {\it Singular Perturbations of Differential Operators.
         Solvable Schr\"{o}dinger Type Operators
    (London Mathematical Society Lecture Note Series vol 271)}
    (Cambridge: Cambridge University Press) p~429


\bibitem{AlbeverioNizhnik:2000}
    Albeverio S and Nizhnik L 2000 
    \emph{Ukrainian Math. J.} \textbf{52} no.~5, 664

\bibitem{Antonevich:1999}
    Antonevich A 1999
    \emph{Nonlinear Phenom. Complex Syst.} \textbf{2} no.~4, 61

\bibitem{BerezinFadeev} Berezin F A and Faddeev L D 1961  {\it Sov. Math. Dokl.} {\bf2} 372

\bibitem{BrascheFigariTeta}
    Brasche J F, Figari R, and Teta A 1998
    \emph{Potential Anal.}  \textbf{8}  no.~2, 163

\bibitem{BrascheNizhnik:2002} Brasche J and Nizhnik L 2002
    \emph{Methods Funct. Anal. Topology}  \textbf{8} no.~3, 13

\bibitem{CacciapuotiExner:2007}
    Cacciapuoti C and Exner P 2007
    \emph{J. Phys. A: Math. Theor.} \textbf{40} no.~26, F511 

\bibitem{ChristianZolotarIermak03}
    Christiansen P L, Arnbak H C, Zolotaryuk A V, Ermakov V N and
    Gaididei Y B 2003 \JPA {\bf 36} 7589

\bibitem{CurgusLangerJDE1989}
     \'{C}urgus B and Langer H   1989 J. Diff. Eq. \textbf{79}, no.~1, 31

\bibitem{Exner:1995} Exner P 1995 {\it J. Math. Phys.} {\bf36} 4561

\bibitem{ExnerNeidhardtZagrebnov}
    Exner P, Neidhardt H, and Zagrebnov V A 2001
    \emph{Comm. Math. Phys.}  \textbf{224} no.~3, 593

\bibitem{GolovatyManko1} Golovaty Yu  and Man'ko S  2009 {\it Ukr. Math. Bulletin} {\bf6} no.~2, 173 (arXiv:0909.1034v1 [math.SP])

\bibitem{GesHolJPA}
    Gesztesy F and Holden H 1987 \JPA  {\bf20} 5157

\bibitem{Friedman:1972}
    Friedman C 1972 \emph{J. Funct. Anal.} \textbf{10} no.~3, 346

\bibitem{Koshmanenko:1999}
    Koshmanenko V 1999 \emph{Singular quadratic forms in perturbation theory.}
    (\emph{Mathematics and its Applications} \textbf{474}) (Dordrecht: Kluwer Academic Publishers) p~308

\bibitem{KostenkoMalamud:2009}
    Kostenko A and Malamud M \rd{2010 \emph{J.~Diff. Equat} \textbf{249} no.~2, 253}

\bibitem{KronigPenney} Kronig R de L and Penney W G 1931
    {\it Proc. Roy. Soc. (London)} {\bf 130A} 499

\bibitem{Kurasov:2003}
    Kurasov P  2003
    {\it Integr. Eq. Oper. Theory} \textbf{45} 437

\bibitem{Manko:2009}
    Man'ko S 2010
    {\it Visnyk of the Lviv University. Mech. and Math.} \textbf{71} 150.

\bibitem{NizhFAA2003}
    Nizhnik L P 2003 {\it Funct. Anal. and Appl.} {\bf37} 85

\bibitem{NizhFAA2006}
    Nizhnik L P 2006 {\it Funct. Anal. and Appl.} {\bf40} 74

\bibitem{Seba:1985}
    \v{S}eba P 1985
    \emph{Lett. Math. Phys.} \textbf{10} no.~1, 21

\bibitem{SebRMP}
    \v{S}eba P 1986 {\it Rep. Math. Phys.} {\bf24} 111

\bibitem{ShkalikovSavchukTMMO2003}
    Savchuk  A M and Shkalikov A A 2003 {\it Trans. Mosc. Math. Soc. } {\bf64} 159

\bibitem{ToyamaNogami}
    Toyama F and Nogami Y 2007 \emph{J. Phys. A: Math. Theor.} {\bf40} F685

\bibitem{ZolotarChristianIermak06}
    Zolotaryuk A V, Christiansen P L, and Iermakova S V   2006 \JPA {\bf 39} 9329

\bibitem{ZolotarChristianIermak07}
    Zolotaryuk A V, Christiansen P L, and Iermakova S V 2007 \emph{J. Phys. A: Math. Theor.} \textbf{40} 5443

\bibitem{Zolotaryuk08}
    Zolotaryuk A V 2008 {\it Adv. Sci. Lett.} {\bf1} 187

\bibitem{Zolotaryuk09}
    Zolotaryuk A V 2010 \emph{J. Phys. A: Math. Theor.} \textbf{43} 105302

\bibitem{Zolotaryuk10}
    Zolotaryuk A V 2010 \emph{Physics Letters A} \rd{\textbf{374} 1636}


\end{thebibliography}
\end{document}